\documentclass[aos,MSNbibl,nameyear,dvips]{arximspdf}
\usepackage{graphicx}

\doi{10.1214/12-AOS1068} 
\volume{41}
\issue{2}
\pubyear{2013}
\firstpage{722}
\lastpage{750}

\makeatletter

\newcommand{\eqref}[1]{(\ref{#1})}
\def\cal{\mathcal}

\newtheorem{theorem}{Theorem}
\newtheorem{corollary}{Corollary}
\newtheorem{proposition}{Proposition}
\newtheorem{lemma}{Lemma}
\newproclaim{example}{Example}
\newproclaim{remark}{Remark}
\newproclaim{definition}{Definition}
\newcommand{\ep}{\varepsilon}
\newcommand{\argmin}{\mathop{\arg\min}}
\def\sf{{\cal F}}
\def\si{{\cal I}}

\def\sf{{\cal F}}
\def\si{{\cal I}}
\def\pr{\mathsf{P}} 

\makeatother

\begin{document}
\begin{frontmatter}

\title{Adaptive confidence intervals for regression functions under
shape constraints}
\runtitle{Adaptive confidence intervals}

\begin{aug}
\author{\fnms{T. Tony} \snm{Cai}\thanksref{t2}\ead[label=e1]{tcai@wharton.upenn.edu}\ead[label=u1,url]{http://www-stat.wharton.upenn.edu/\textasciitilde tcai}},
\author{\fnms{Mark G.} \snm{Low}\corref{}\ead[label=e2]{lowm@wharton.upenn.edu}\ead[label=u2,url]{http://www-stat.wharton.upenn.edu/\textasciitilde lowm}}
\and
\author{\fnms{Yin} \snm{Xia}\thanksref{t2}\ead[label=e3]{xiayin@wharton.upenn.edu}\ead[label=u3,url]{http://www-stat.wharton.upenn.edu/\textasciitilde xiayin}}
\thankstext{t2}{Supported in part by NSF FRG Grant DMS-08-54973.}

\runauthor{T. T. Cai, M. G. Low and Y. Xia}
\affiliation{University of Pennsylvania}
\address{Department of Statistics\\
The Wharton School\\
University of Pennsylvania\\
Philadelphia, Pennsylvania 19104\\
USA\\
\printead{e1}\\
\phantom{E-mail:\ }\printead*{e2}\\
\phantom{E-mail:\ }\printead*{e3}\\
\printead{u1}\\
\phantom{URL:\ }\printead*{u2}\\
\phantom{URL:\ }\printead*{u3}}
\end{aug}

\received{\smonth{5} \syear{2012}}
\revised{\smonth{11} \syear{2012}}

%
\begin{abstract}
Adaptive confidence intervals for regression functions are constructed
under shape constraints of monotonicity and convexity.
A~natural benchmark is established for the minimum expected length of
confidence intervals at a given function in terms of an analytic
quantity, the local modulus of continuity.
This bound depends not only on the function but also the assumed
function class.
These benchmarks show that the constructed confidence intervals have
near minimum expected length for each individual function, while
maintaining a given coverage probability for functions within the class.
Such adaptivity is much stronger than adaptive minimaxity over a
collection of large parameter spaces.
\end{abstract}

%

\end{frontmatter}

\section{Introduction}
\label{intro.sec}

The construction of useful confidence sets is one of the more challenging
problems in nonparametric function estimation.
There are two main interrelated issues which need to
be considered together, coverage probability and
the expected size of the confidence set.
For a fixed parameter space it is often possible
to construct confidence sets which have guaranteed coverage probability
over the parameter
space while controlling the maximum expected size.
However such minimax statements
are
often thought to be too conservative, and
a more natural goal is to have the expected size of the confidence
set reflect in some sense the difficulty of estimating
the particular underlying function.

These issues are well illustrated by considering
confidence intervals for the value of a function at a fixed point.
Let $Y$ be an observation from the
white noise model
%
\begin{equation}
dY(t) = f(t) \,dt + n^{-{1}/{2}} \,dW(t),\qquad -\tfrac{1}{2}\le t \le
\tfrac{1}{2}, \label{w.model}
\end{equation}
where $W(t)$ is standard Brownian motion and $f $ belongs to some
parameter space~${\sf}$.
Suppose\vspace*{1pt} that we wish to construct a confidence interval for $f$ at some
point $t_0\in(-\frac{1}{2}, \frac{1}{2})$.
Let $\mathrm{CI}$ be a confidence interval for $f(t_0)$ based on observing the
process $Y$, and let $L(\mathrm{CI})$ denote the length of the confidence interval.
The minimax point of view can then be expressed by the following:
subject to the constraint on the coverage probability $\inf_{f \in\sf
}P( f(t_0) \in \mathrm{CI} ) \ge1 - \alpha$,
minimize the maximum expected length $\sup_{f \in\sf} E_f( L(\mathrm{CI}))$.

As an example it is common to consider the Lipschitz classes
\[
\Lambda(\beta,M) = \bigl\{ f\dvtx\bigl |f(y) - f(x) \bigr| \le M |y-x|^{\beta}
\mbox{ for $x, y\in\bigl[-\tfrac{1}{2}, \tfrac{1}{2}\bigr]$} \bigr\},\qquad
\mbox{if $0<\beta\le1$}
\]
and for $\beta> 1$
\[
\Lambda(\beta, M) = \bigl\{ f\dvtx \bigl|f^{(\lfloor\beta\rfloor)}(x) - f^{(\lfloor\beta\rfloor)}(y)\bigr| \leq
M |x-y|^{\beta'} \mbox{ for $x, y\in\bigl[-\tfrac{1}{2},
\tfrac
{1}{2}\bigr]$} \bigr\},
\]
where $\lfloor\beta\rfloor$ is the largest integer less than $\beta$
and $\beta' = \beta- \lfloor\beta\rfloor$.
For these classes it easily follows from results of \citet{Donoho}, \citet{Low97} and \citet{Evans} that the minimax expected
length of confidence intervals, which have guaranteed coverage of $1-
\alpha$ over $\Lambda(\beta,M)$, is of order $M^{1 /(1 + 2 \beta)}n^{-{\beta/(1 + 2\beta)}}$.


It should, however, be stressed that confidence intervals which achieve
such an expected length rely on the knowledge of the particular
smoothness parameters $\beta$ and $M$, which are not known in most
applications.
Unfortunately, \citet{Low97} and \citet{CL04} have shown that the
natural goal of constructing an adaptive confidence interval which has
a given coverage probability and
has expected length that is simultaneously close to these minimax
expected lengths
for a range of smoothness parameters is not in general attainable.
More specifically suppose that
a confidence interval has guaranteed coverage probability of $1- \alpha
$ over $\Lambda(\beta,M)$.
Then for any $f$ in the interior of $\Lambda(\beta, M)$
the expected length for this $f$ must also be of order
$n^{-{\beta/(1 + 2\beta)}}$.
In other words the minimax rate describes the actual rate for all
functions in the class
other than those on the boundary of the set.
For example, in the case that
a confidence interval has guaranteed coverage probability of $1 -
\alpha
$ over the Lipschitz class $\Lambda(1,M)$, then even if the underlying
function has two derivatives, and the first derivative smaller than~$M$,
the confidence interval
for $f(x)$ \textit{must} still have expected length of order $n^{-1/3}$
even though one would hope
that an adaptive confidence interval would have a much shorter length
of order $n^{-2/5}$.

Despite these very negative results there are some settings where some
degree of adaptation has been shown to be possible.
In particular under certain shape constraints \citet{HenSta95}
constructed confidence bands which have a guaranteed coverage
probability of at least $1 - \alpha$ over the collection of all
monotone densities and which have maximum expected length of order
$(\frac{\log n}{n})^{ \beta/(2 \beta+1)}$ for those monotone
densities which are in $\Lambda(\beta,M)$ for a particular choice of
$\beta$ where $0 < \beta\le1$.
This construction relies on the selection of a tuning parameter and is
thus not adaptive.
D\"{u}mbgen (\citeyear{Dum03}), however, does provide adaptive confidence bands
with optimal rates for both isotonic and convex functions under
supremum norm loss on arbitrary compact subintervals.
These results are, however, still framed in terms of the maximum length
over particular large parameter spaces,
and the existence of such intervals raises
the question of exactly how much adaption is possible.
It is this question that is the focus of the present paper.

Rather than considering the maximum expected length over large
collections of functions,
we study the problem of adaptation to each and every function in the
parameter space.
We examine this problem in detail
for two commonly used collections of functions that have shape
constraints, namely the collection of convex functions and the
collection of monotone functions.
We focus on these parameter spaces as it is for such shape constrained
problems for which there is some hope
for adaptation.
Within this context we consider the problem of constructing a
confidence interval for the value of a function at a fixed point under
both the white noise with drift model given in (\ref{w.model}) as well
as a nonparametric regression model.
We show that within the class of convex functions and the class of
monotone functions, it is indeed possible to \textit{adapt to each
individual function}, and not just to the minimax expected length over
different parameter spaces in a collection.
The notion of adaptivity to a single function is also discussed in
Lepski, Mammen and Spokoiny (\citeyear{LepMamSpo97}) and \citet{LepSpo97} for
the related point estimation problem but in these contexts a
logarithmic penalty of the noise level must be paid, and thus
the notion of adaptivity is somewhat different.\looseness=1

This result is achieved in two steps. First we study the problem of
minimizing the expected length of
a confidence interval, assuming that the data is generated from a
particular function $f$ in the parameter space,
subject to the constraint that the confidence interval has guaranteed
coverage probability over the entire parameter space.
The solution to this problem gives a benchmark for the expected length
which depends on the function $f$ considered.
It gives a bound on the expected length of any adaptive interval because
if the expected length is smaller than this bound for \textit{any}
particular function, the confidence interval cannot have the desired
coverage probability.
In applications it is more useful to express the benchmark in terms of
a local modulus of continuity, an analytic quantity that can be easily
calculated for individual functions.
In situations where adaptation is not possible, this local modulus of
continuity does not vary significantly from function to function. Such
is the case
in the settings considered in \citet{Low97}. However, in the context of
convex or monotone functions,
the resulting benchmark does vary significantly, and this opens up the
possibility for adaptation in those settings.

Our second step is to actually construct adaptive confidence intervals.
This is done separately for monotone functions and convex functions,
with similar results. For example,
an adaptive confidence interval is constructed which is shown to have
expected length uniformly within an absolute constant factor of the
benchmark for every convex function, while maintaining coverage
probability over the collection of all convex functions.
In other words, this confidence interval has smallest expected length,
up to a universal constant factor, for each and every convex function
within the class of all confidence intervals which guarantee a $1-
\alpha$ coverage probability over all convex functions.
A similar result is established for a confidence interval designed for
monotone functions.


The rest of the paper is organized as follows. In Section \ref
{lowerbound} the benchmark for the expected length at each monotone
function or each convex function is established under the constraint
that the interval has a given level of coverage probability over the
collection of monotone functions or the collection of convex functions.
Section~\ref{procedure.sec} constructs data driven confidence intervals
for both monotone functions and convex functions and shows that these
confidence intervals maintain coverage probability and have expected
length within an absolute constant factor of the benchmark given in
Section~\ref{lowerbound} for each monotone function and convex
function. Section~\ref{regression.sec} considers the nonparametric
regression model, and Section~\ref{discussion.sec} discusses
connections of our results with other work in the literature. Proofs
are given in Section~\ref{proof.sec}.

\section{Benchmark and lower bound on expected length}
\label{lowerbound}

As mentioned in the \hyperref[intro.sec]{Introduction}, the focus in this paper is the
construction of confidence intervals which have expected length that
adapts to the unknown
function. The evaluation of these procedures depends on lower bounds
which are given here
in terms of a local modulus of continuity first
introduced by \citet{CL11} in the context of point estimation of
convex functions under mean squared error loss. These lower bounds
provide a natural benchmark for our problems.

\subsection{Benchmark and lower bound}

We focus in this paper on estimating the function $f$ at $0$ since
estimation at other points away from the boundary is similar.
For a given function class $\sf$, write $\si_{\alpha}(\sf)$ for the
collection of all confidence intervals which cover $f(0)$ with
guaranteed coverage probability of $1-\alpha$ for all functions in
$\sf$. For a given confidence interval $\mathrm{CI}$, denote by
$L(\mathrm{CI})$ the length
of $\mathrm{CI}$ and $L(\mathrm{CI}, f) = E_f(L(\mathrm{CI}))$ the expected length of $\mathrm{CI}$ at a
given function~$f$. The minimum expected length at $f$ of all
confidence intervals with guaranteed coverage probability of $1-\alpha$
over $\sf$ is then given by\looseness=-1
%
\begin{equation}
\label{constrained.def}
L_{\alpha}^*(f, \sf) = \inf_{\mathrm{CI} \in\si_{\alpha}(\sf)}L(\mathrm{CI},f).
\end{equation}\looseness=0
A natural goal is to construct a confidence interval with expected
length close to the minimum $L_{\alpha}^*(f, \sf) $ for every $f\in
\sf
$ while maintaining the\vadjust{\goodbreak} coverage probability over $\sf$.
However although $L_{\alpha}^*(f, \sf) $ is a natural benchmark for the
expected length of confidence intervals, it is not easy to evaluate
exactly. Instead as a first step toward our goal, we provide a lower
bound for the benchmark $L_{\alpha}^*(f, \sf) $
in terms of a local modulus of continuity $\omega(\varepsilon,f, \sf)$
introduced by \citet{CL11}. The local modulus is a quantity that
is more easily computable and techniques for its analysis are similar
to those given in \citet{GR3} and \citet{Donoho} where a
global modulus of continuity was introduced in the study of minimax
theory for estimating linear functionals. See the examples in Section~\ref{example.sec}.

For a parameter space $\cal{F}$ and function $f \in\cal{F}$, the local
modulus of continuity is defined by
%
\begin{equation}
\omega(\ep, f, \sf) = \sup\bigl\{ \bigl|g(0) - f(0)\bigr| \dvtx \|g - f
\|_2 \le\ep, g \in\sf\bigr\},
\end{equation}
where $\| \cdot\|_2$ is the $L_2(-\frac{1}{2}, \frac{1}{2})$
function norm.
The following theorem gives a lower bound for the minimum expected
length $L_{\alpha}^*(f, \sf) $ in terms of the local modulus of
continuity $\omega(\varepsilon,f, \sf)$.
In this theorem and throughout the paper we write $\Phi$ for the
cumulative distribution function and $\phi$ for the density function of
a standard normal density and set
$z_{\alpha} = \Phi^{-1}(1 - \alpha)$.
%
\begin{theorem}
\label{EL.bound.thm}
Suppose $\sf$ is a nonempty convex set. Let $0 < \alpha< {1 \over2}$
and $f \in\sf$.
Then for confidence intervals based on (\ref{w.model}),
%
\begin{equation}
\label{EL.bound2} L^*_\alpha(f, \sf) \ge\biggl( 1 - \frac{1}{\sqrt{2\pi} z_{\alpha}} +
\frac
{\phi(z_{\alpha})}{z_{\alpha}} - \alpha\biggr) \omega\biggl(\frac{z_{\alpha
}}{\sqrt
n}, f, \sf\biggr).
\end{equation}
%
In particular,
%
\begin{equation}
L^*_\alpha(f, \sf) \ge\biggl( 1 - \frac{1}{\sqrt{2\pi} z_{\alpha}}\biggr) \omega
\biggl(\frac{z_{\alpha}}{\sqrt n}, f, \sf\biggr).
\end{equation}
%
\end{theorem}

The lower bounds given in Theorem~\ref{EL.bound.thm} can be viewed as
benchmarks for the evaluation of the expected length of confidence
intervals when the true function is~$f$
for confidence intervals which have guaranteed coverage probability
over all of $\sf$.
The bound depends on the underlying true function $f$ as well as the
parameter space $\sf$.

The bounds from Theorem~\ref{EL.bound.thm} are general. In some
settings they can be used to rule out the possibility of adaptation,
whereas in other settings they provide bounds on how much adaptation is
possible.
In particular the result ruling out adaptation over Lipschitz classes
mentioned in the \hyperref[intro.sec]{Introduction} easily follows from this theorem.
For example, consider the Lipschitz class $\Lambda(\beta, M)$ and
suppose that $f$ is in the interior of $\Lambda(\beta,M)$.
Straightforward calculations similar to those given in Section \ref
{example.sec} show that
%
\begin{equation}
\label{Lip.modulus} \omega\bigl(\ep, f, \Lambda(\beta, M)\bigr) \sim C
\ep^{2\beta/{(2 \beta+1)}}.
\end{equation}


Now consider two Lipschitz classes $\Lambda(\beta_1,M_1)$ and
$\Lambda
(\beta_2,M_2)$ with $\beta_1> \beta_2$.
A~fully adaptive confidence\vadjust{\goodbreak} interval in this setting would have
guaranteed coverage of $1-\alpha$ over $\Lambda(\beta_1,M_1)\cup
\Lambda
(\beta_2,M_2)$ and maximum expected length over $\Lambda(\beta_i,M_i)$
of order $n^{\beta_i/(2\beta_i + 1)}$ for $i=1$ and $2$. However, it
follows from Theorem~\ref{EL.bound.thm} and \eqref{Lip.modulus} that
for all confidence intervals with coverage probability of $1-\alpha$
over $\Lambda(\beta_2,M_2)$, for every $f\in\Lambda(\beta_2,M')$
with $M'<M_2$,
\[
L^*_\alpha\bigl(f, \Lambda(\beta_2,M_2)\bigr)
\ge C(\alpha) n^{-{\beta_2 /(2 \beta_2 +1)}}
\]
for some constant $C(\alpha)$ not depending on $f$. In particular this
holds for all $f\in\Lambda(\beta_1,M_1)\cap\Lambda(\beta_2,M')$
and hence
\begin{eqnarray*}
\sup_{f\in\Lambda(\beta_1,M_1)}\inf_{\mathrm{CI}\in\si_\alpha(\Lambda
(\beta
_1,M_1)\cup\Lambda(\beta_2,M_2))} L(\mathrm{CI}, f) &\ge& C(
\alpha)n^{-{\beta_2 /(2 \beta_2 +1)}}\\
&\gg& n^{-{\beta_1 /(2 \beta_1 +1)}}.
\end{eqnarray*}
Therefore it is not possible to have confidence intervals with adaptive
expected length over two Lipschitz classes with different smoothness parameters.

In the present paper
Theorem~\ref{EL.bound.thm}
will be used to provide benchmarks in the setting
of shape constraints.
Denote by $F_m$ and $F_c$, respectively, the collection of all
monotonically nondecreasing functions and the collection of all convex
functions on $[-\frac{1}{2}, \frac{1}{2}]$.
We shall now show that
in these cases
the modulus and the associated lower bounds vary significantly from
function to function.


\subsection{Examples of bounds for monotone functions and convex functions}
\label{example.sec}

We now turn to the application of the lower bound given in Theorem \ref
{EL.bound.thm}
in the case of monotone functions and convex functions.
Here we shall evaluate the lower bound for
four particular families of functions yielding different rates at which
the expected length decreases to zero as
the noise level decreases in contrast to the situation just described
where the parameter space did not
have an order constraint.
Two of the functions will be both monotonically nondecreasing and
convex. In this case
the lower bound can also be quite different depending on whether we
assume the knowledge that
$f$ is convex or monotonically nondecreasing.

The key quantity that is needed in any application of Theorem \ref
{EL.bound.thm} is the local modulus.
We follow the same approach as given in \citet{Donoho} where a global
modulus of continuity is considered for minimax estimation.
In each case, for a given function $f$, we first minimize the $L_2$
norm between a function $g \in{\cal F}$ and the function $f$ subject
to the constraint that $|g(0) - f(0)| = a$ for some given value $a>0$.
From here it is easy to invert and thus maximize $|g(0) - f(0)|$ given
a constraint on the $L_2$ norm between $f$ and $g$.

\begin{example}
\label{linear.ex}
As a first example consider the linear function $f_k(t)=kt$ where $k
\ge0$ is a constant.
This function is both monotonically nondecreasing and convex.\vadjust{\goodbreak}

First consider the collection of monotonically nondecreasing functions
$F_m$. We shall treat separately the case $k>0$ and the case $k=0$. For
the moment we shall take $k >0$.
Suppose that $0 < a \le\frac{k}{2}$.
In this case $f_k\in F_m$ and a
function $g$ that minimizes $\|g-f_k\|_2$ subject to the constraint
that $|g(0) - f_k(0)| = a$ is given by
$g(t) = f_k(t)$ if $t <0$, $g(t) = a$ if $0 \le t \le b$, and $g(t) =
f_k(t)$ if $t >b$,
where $b$ satisfies $f_k(b) = a$. The assumption that $a \le\frac
{k}{2}$ guarantees
$b \le\frac{1}{2}$.
We then have $\|g - f_k\|_2 = a^{{3}/{2}}/(3k)^{{1}/{2}} $.
It follows that if $\varepsilon^2\leq\frac{1}{24}k^2$
\[
\omega(\varepsilon,f_k, F_m) = (3k)^{{1}/{3}}
\varepsilon^{{2}/{3}}
\]
and consequently for $n\geq\frac{24z_{\alpha}^2}{k^2}$, if $k >0$
\[
L^*_\alpha(f_k, F_m) \ge\biggl( 1 -
\frac{1}{\sqrt{2\pi} z_{\alpha}}\biggr) (3k)^{{1}/{3}}z^{2/3}_{\alpha}n^{-{1/3}}.
\]

In the case that $k=0$
a function $g$ that minimizes $\|g-f_0\|_2$ subject to the constraint
that $|g(0) - f_0(0)| = a$ is given by
$g(t) = f_0(t)$ if $t <0$, $g(t) = a$ if $0 \le t \le\frac{1}{2}$.
In this case it is easy to check that $\|g - f_0\|_2 = \frac{1}{\sqrt
2} \varepsilon$ and hence
\[
\omega(\varepsilon,f_0, F_m) = \sqrt{2} \varepsilon
\]
and hence
\[
L^*_\alpha(f_0, F_m) \ge\biggl( 1 -
\frac{1}{\sqrt{2\pi} z_{\alpha
}}\biggr)\sqrt {2} z_{\alpha}n^{-{1/2}}.
\]

We now consider the bound for the length of the confidence interval
for $f_k$ belonging to the collection of convex functions.
In this case we do not need to treat the cases $k>0$ and $k=0$ separately.
The function $g$ that minimize $\|g-f_k\|_2$ subject to the constraint
that
$g$ is convex and $|g(0)- f_k(0)|=a$ is given by $g(t)=(k+3a)t-a$ if
$t\geq0$ and $g(t)=(k-3a)t-a$ if $t<0$. In this case $\|g-f\|_2=\frac
{1}{2}a$. It then immediately follows that
\[
\omega(\varepsilon, f_k, F_c) = 2 \varepsilon
\]
and so
\[
L^*_\alpha(f_k, F_c) \ge\biggl( 1 -
\frac{1}{\sqrt{2\pi} z_{\alpha}}\biggr) 2z_{\alpha}n^{-{1/2}}.
\]

It is important to note that for $k>0$ the minimum expected lengths
$L^*_\alpha(f_k, F_m) $ and $L^*_\alpha(f_k, F_c) $ are different, one
of order $n^{-{1/3}}$ and another of order $n^{-{1/2}}$,
although the function $f_k$ is the same.
It is also interesting to note that the expected length of the
confidence for monotone functions is an increasing function of $k$
whereas the expected length of the confidence for convex functions does
not depend on $k$.
Since we shall show that these bounds are achievable within a constant factor
it follows that the minimum expected length of the confidence interval
when $f_k$ is the true function
depends strongly on whether we specify that the underlying collection
of functions is convex or monotone. Plots illustrating shapes of
functions $f_k$ and a least favorable function $g$ are shown as below
in Figure~\ref{fig:ex1}.\looseness=1

\begin{figure}

\includegraphics{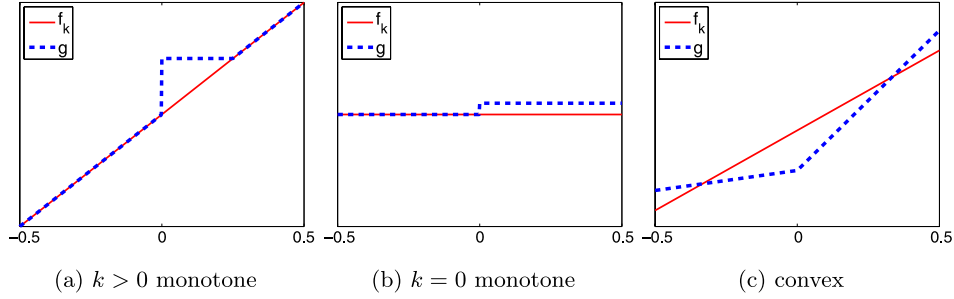}

\caption{Plots of $f_k$ and a least favorable function $g$ in Example
\protect\ref{linear.ex} with the constraints $|g(0)- f_k(0)|=a$.}\label{fig:ex1}
\end{figure}

%
%

\end{example}

\begin{example}\label{ex2}
As a second example\vspace*{1pt} which is also both monotonically nondecreasing and convex
consider the function $f(t)= k_1t + k_2t^r I\ (0<t\le\frac{1}{2})$
where $r\ge
1$ and $k_1 \ge0$ and $k_2 > 0$ are constants.

We consider the cases $r=1$ and $r>1$ separately. When $r=1$ the
function is piecewise linear with the change
of slope at $0$.
In this case suppose $0<a\leq\frac{k_1+k_2}{2}$. A monotonically
nondecreasing function $g\in F_m$ that minimize $\|g-f\|_2$ subject to
the constraint that $|g(0)-f(0)|=a$ is given by $g(t)=f(t)$ if $t<0$,
$g(t)=a$ if $0\leq t\leq b$, and $g(t)=f(t)$ if $t>b$, where $b$
satisfies $f(b)=a$. The constraint $a \leq\frac{k_1 +k_2}{2}$ is to
guarantee that such a $b$ exists with $b \le\frac{1}{2}$.
Then we have $\|g-f\|_2=a^{{3}/{2}}(3(k_1+k_2))^{-{1}/{2}}$,
and it follows that if $\varepsilon^2\leq\frac{1}{24}(k_1+k_2)^2$,
\[
\omega(\varepsilon, f, F_m) = \bigl(3(k_1+k_2)
\bigr)^{{1}/{3}}\varepsilon ^{{2}/{3}}
\]
and consequently for $n\geq\frac{24z_{\alpha}^2}{(k_1+k_2)^2}$,
\[
L^*_\alpha(f, F_m) \ge\biggl( 1 - \frac{1}{\sqrt{2\pi} z_{\alpha}}\biggr)
\bigl(3(k_1+k_2)\bigr)^{{1}/{3}}z^{2/3}_{\alpha}n^{-{1/3}}.
\]

We can also give a lower bound on the expected length for this same
function for confidence intervals which guarantee
coverage over the class of convex functions.
Suppose $0<a\leq\frac{k_2}{4}$. Here we need to find the convex $h$
that minimizes $\|h-f\|_2$ subject to the constraints that
$|h(0)-f(0)|=a$. It is given by $h(t)=f(t)$ if $t\leq-\frac{2a}{k_2}$,
$h(t)=(\frac{k_2}{2}+k_1)t+a$ if $-\frac{2a}{k_2}\leq t\leq\frac
{2a}{k_2}$ and $h(t)=f(t)$ if $t\geq\frac{2a}{k_2}$. Then $\|f-g\|
_2=2a^{{3}/{2}}/(3k_2)^{{1}/{2}}$ and it follows that if
$\varepsilon^2\leq\frac{k_2^2}{48}$,
\[
\omega(\varepsilon,f,F_c) = (3k_2/4)^{{1}/{3}}
\varepsilon^{{2}/{3}}.
\]
Hence, for $n\geq\frac{48z_{\alpha}^2}{k_2^2}$,
\[
L^*_\alpha(f,F_c)\geq\biggl( 1 - \frac{1}{\sqrt{2\pi} z_{\alpha}}\biggr)
(3k_2/4)^{{1}/{3}}z^{2/3}_{\alpha}n^{-{1/3}}.
\]

We now turn to the case where $r>1$.
Suppose $0<a\leq\frac{k_1}{2}+k_2(\frac{1}{2})^r$.
In this case the monotonically nondecreasing function $g$ that
minimizes $\|g-f\|_2$ subject to the constraints that $|g(0)-f(0)|=a$
is given by $g(t)=f(t)$ if $t<0$, $g(t)=a$ if $0\leq t\leq b$ and
$g(t)=f(t)$ if $t>b$, where $b$ satisfies $f(b)=a$. As before the condition
$0<a\leq\frac{k_1}{2}+k_2(\frac{1}{2})^r$ guarantees that $b$ exists
with $b < \frac{1}{2}$.
In this case
$a^{{3}/{2}}(3k_1)^{-{1}/{2}}- ca^s \le\|g-f\|_2 \le
a^{{3}/{2}}(3k_1)^{-{1}/{2}}+ca^s$
for some constant $c>0$ and $s>3/2$. It follows that if $\varepsilon
^2\leq\frac{1}{24}k_1^2+(1+\frac{1}{2r+1}-\frac{2}{r+1})(\frac
{1}{2})^{2r+1}k_2^2+(\frac{1}{2}-\frac{1}{r+1}+\frac{1}{r+2})(\frac
{1}{2})^{r+1}k_1k_2$, then
\[
\omega(\varepsilon, f, F_m) = (3k_1)^{{1}/{3}}
\varepsilon^{{2}/{3}}\bigl(1+o(1)\bigr).
\]
Hence,
%
\[
L^*_\alpha(f, F_m) \ge\biggl( 1 - \frac{1}{\sqrt{2\pi} z_{\alpha}}\biggr)
(3k_1)^{{1}/{3}}z^{2/3}_{\alpha}n^{-{1/3}}
\bigl(1+o(1)\bigr).
\]

For a bound on the expected length of this same function for confidence
intervals with coverage
guaranteed over the collection of convex functions, we suppose $0<a\leq
k_2(\frac{1}{2})^{r+1}$. In this case
the convex function $h$ that minimizes $\|h-f\|_2$ subject to the
constraints that $|h(0)-f(0)|=a$, is given
by $h(t)=kt+a$, $k>k_1$, if $x_0\leq t\leq x_1$ and $h(t)=f(t)$
otherwise, where $(x_0,cx_0)$ and $(x_1,cx_1+x_1^r)$ are the
intersection points of $f(t)$ and the line $kt+a$. Then the function
$h$ with slope $k_0$ that minimize $\|h-f\|_2$ would be the least
favorable function. It follows that, if $\varepsilon^2\leq\frac
{k_2^2}{24}(\frac{1}{2})^{2r}$,
\[
\omega(\varepsilon, f, F_c) = C(r)k_2^{{1}/{(2r+1)}}
\varepsilon^{{2r}/{(2r+1)}}
\]
and consequently for $n\geq\frac{24z_{\alpha}^22^{2r}}{k_2^2}$,
\[
L^*_\alpha(f, \sf_c) \ge\biggl( 1 - \frac{1}{\sqrt{2\pi} z_{\alpha
}}
\biggr)C(r)k_2^{{1}/{(2r+1)}} z_{\alpha}^{{2r}/{(2r+1)}}
n^{-{r/(2r+1)}},
\]
where $C(r)>0$ is a constant depending on $r$ only.
\end{example}

It is interesting to note that in this example the rates of convergence
for $L^*_\alpha(f, \sf_m)$ and $L^*_\alpha(f, \sf_c)$ are the same for
the case $r=1$, and are different when $r>1$. Plots illustrating shapes
of functions $f$ and a least favorable function $g$ are shown as below
in Figure~\ref{fig:ex2}.\vadjust{\goodbreak}

\begin{figure}

\includegraphics{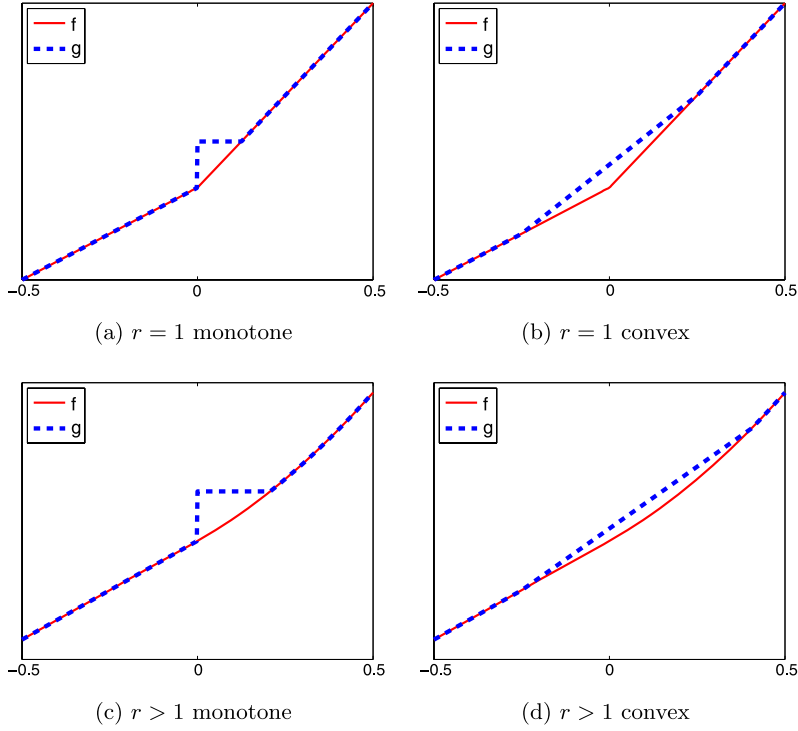}

\caption{Plots of $f$ and a least favorable function $g$ in Example \protect\ref{ex2}
with the constraints $|g(0)- f(0)|=a$.}\label{fig:ex2}
\end{figure}

%

Next we consider a function which is monotonically nondecreasing but
not convex.

\begin{example}\label{ex3}
Let $f(t)=kt^r$ for some constant $k > 0$ and $r=2l+1$ or $r=\frac
{1}{2l+1}$ for $l=0,1,2,\ldots.$
Suppose that $a<(\frac{1}{2})^rk$. In this case a function $g$ that
minimizes $\Vert g-f\Vert _2$ subject to the constraint
that $|g(0) - f(0)| = a$ is given by
$g(t) = f(t)$ if $t <0$, $g(t) = a$ if $0 \le t \le b$ and $g(t) =
f(t)$ if $t >b$,
where $b$ satisfies $f(b) = a$. As before the condition
$a<(\frac{1}{2})^rk$ guarantees that $b$ exists with $b < \frac{1}{2}$.
Then $\Vert g - f\Vert _2 = 
a^{1+({1}/{(2r)})}/k^{{1}/{(2r)}}(2r^2/(r+1)(2r+1))^{{1}/{2}}$, and
it follows that when $\varepsilon^2\leq(\frac{1}{2})^{2r+1}k^2\frac
{2r^2}{(r+1)(2r+1)}$,
\[
\omega(\varepsilon,f, F_m) = \biggl(\frac{(r+1)(2r+1)k}{2r^2}
\biggr)^{{r}/{(2r+1)}} \varepsilon^{{2r}/{(2r+1)}}.
\]
Hence for $n\geq\frac{2^{2r+1}(r+1)(2r+1)z_{\alpha}^2}{2r^2k^2}$,
\[
L^*_\alpha(f, \sf_c) \ge\biggl( 1 - \frac{1}{\sqrt{2\pi} z_{\alpha}}
\biggr) \biggl(\frac
{(r+1)(2r+1)k}{2r^2}\biggr)^{{r}/({2r+1)}} z_{\alpha}^{2r /(2r+1)}
n^{-{r/(2r+1)}},
\]
and once again it is clear that the rate at which the expected length
decreases to zero depends strongly on the value of $r$.
\end{example}

As a final example we consider a function which is convex but not
monotonically nondecreasing.

\begin{example}\label{ex4}
Let $f(t)=t^2$.
Suppose that $a<1/2$. In this case the function $g$ that minimizes
$\Vert g-f\Vert _2$ subject to the constraint
that $|g(0) - f(0)| = a$ is given by
$g(t) = -3\sqrt{a/2}t-a$ if $-\sqrt{2a}\leq t \leq0$, $g(t) = 3\sqrt {a/2}t-a$ if $0 \le t \le\sqrt{2a}$ and $g(t) = f(t)$ otherwise.
Then $\Vert g - f\Vert _2 = 2^{5/4}/\sqrt{15}a^{5/4}$ and
it follows that when $\varepsilon^2\leq1/\sqrt{15}$,
\[
\omega(\varepsilon,f, F_m) = \frac{15^{2/5}}{2}\varepsilon^{{4}/{5}}.
\]
Hence for $n\geq\sqrt{15}z_{\alpha}^2$,
\[
L^*_\alpha(f, \sf_c) \ge\biggl( 1 - \frac{1}{\sqrt{2\pi} z_{\alpha}}
\biggr) \frac
{15^{2/5}}{2}z_{\alpha}^{{4}/{5}}n^{-{2}/{5}}.
\]
A similar minimization problem is solved in D\"{u}mbgen (\citeyear{Dum03}).
\end{example}

\begin{figure}

\includegraphics{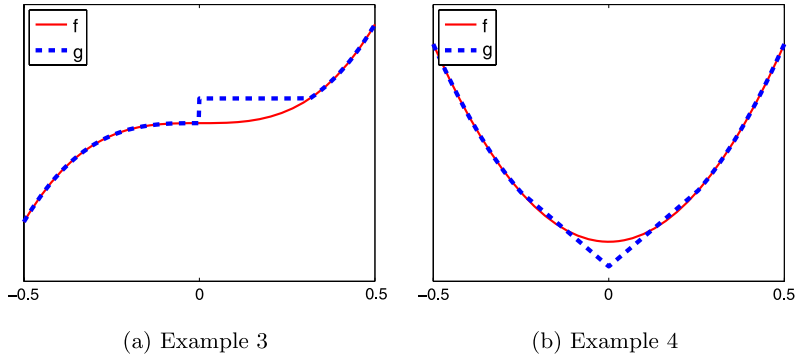}

\caption{Plots of $f$ and a least favorable function $g$ in Examples \protect\ref{ex3}
and \protect\ref{ex4} with the constraints $|g(0)- f(0)|=a$.}\label{fig:ex3-4}
\end{figure}

Plots illustrating shapes of functions $f$ and a least favorable
function $g$ for both Examples~\ref{ex3} and~\ref{ex4} are shown in Figure
\ref{fig:ex3-4}.

%


\section{Confidence procedures}
\label{procedure.sec}

In this section we both construct and give an analysis of adaptive
confidence intervals for monotone functions and convex functions. The
procedures are easily implementable. We consider the class of
monotonically nondecreasing functions and the class of convex
functions. Concave functions and monotonically nonincreasing functions
can be handled similarly.

\subsection{Construction}

The construction is split into two steps. In the first step a countable
collection of confidence intervals is created each of which has
guaranteed coverage probability. These intervals are based on a
collection of pairs of linear estimators. For each interval one of the
estimators has nonnegative bias and the other nonpositive bias.
The one-sided control of the bias of these estimators is a key special
feature in these problems
and an important part of what makes it possible to adapt to every
individual function.
Moreover for each function $f$ this collection has at least one
interval with expected length within a constant factor of the local
modulus bound given in Theorem~\ref{EL.bound.thm}.
The second step is to select from this collection a particular interval.



In the case of monotonically nondecreasing functions we take
for each $j\geq2$, pairs of estimators $\delta_j^R=2^{j}(Y(2^{-j}) -
Y(0))$ and $\delta_j^L=2^{j}(Y(0) - Y(-2^{-j}))$.
Then for estimating $f(0)$ it is easy to check
that $\delta_j^R$ has nonnegative and monotonically nonincreasing
biases while $\delta_j^L$ have nonpositive and monotonically
nondecreasing biases. 
The one-sided control of the biases of these estimators over the class
of all monotonically nondecreasing functions
easily allows for the construction of a confidence interval.
For that we shall need the standard deviation of $\delta_j^R$ and
$\delta_j^L$.
In order to give a unified treatment in both the monotone and convex
case it is useful
to establish a common notation. Here we shall set
$\sigma^2_j=\frac{2^{j-1}}{n}$.
It is then easy to check that both $\delta_j^R$ and $\delta_j^L$ have a
standard deviation of $ \sqrt2 \sigma_j$.
It
is then also easy to see that
for each $j \ge2$, the confidence interval
$\mathrm{CI}^m_j(\alpha)$
given by
%
\begin{equation}
\label{mCIj} \mathrm{CI}^m_j(\alpha)=\bigl[\delta_j^L-z_{\alpha/2}
\sqrt{2}\sigma_j,\delta _j^R+z_{\alpha/2}
\sqrt{2}\sigma_j\bigr]
\end{equation}
has guaranteed coverage of $1-\alpha$.
We should, however, note that in (\ref{mCIj})
the left endpoint of the interval may be larger than the right endpoint
in which case we
adopt the convention that the confidence interval is just the empty
set. The length of this confidence interval is then $\max(\delta
_j^R-\delta_j^L+2\sqrt{2}z_{\alpha/2}\sigma_j,0)$.

In the case of convex functions
for $j \ge1$, let $\delta_j = 2^{j-1}(Y(2^{-j}) - Y(-2^{-j}))$
and let $\tilde\delta_j = 2\delta_{j+1} - \delta_j$.

The following lemma shows that for convex functions $\delta_j$ have
nonnegative and monotonically
nonincreasing biases and that $\tilde\delta_j$
have nonpositive and monotonically nondecreasing biases.
%
\begin{lemma}\label{con.lem}
For any convex function $f$,
%
\begin{eqnarray}
0&\leq&\operatorname{Bias}(\delta_{j+1})\leq\tfrac{1}{2}
\operatorname {Bias}(\delta _j),\label{con.bias.eq}
\\
E\delta_j-3E\delta_{j+1}+2E\delta_{j+2}&\geq&0
\label{con.eq}.
\end{eqnarray}
\end{lemma}

It is also easy to check that
the standard deviation of $\delta_j$ is equal to $\sigma_j$
where
$\sigma^2_j=\frac{2^{j-1}}{n}$
and that\vadjust{\goodbreak}
$2\delta_{j+1} - \delta_j$ has a standard deviation of ${\sqrt5}
\sigma_j$.
It then follows from the signs of the biases of
$\delta_{j+1}$ and $2\delta_{j+1} - \delta_{j}$ that
for any given $j$,
%
\begin{equation}
\label{CIsj} \mathrm{CI}_j^c(\alpha) = [2 \delta_{j+1}
- \delta_j - z_{\alpha/2} {\sqrt{5}} \sigma_j,
\delta_{j+1} + z_{\alpha/2} \sigma_{j+1}]
\end{equation}
gives a confidence interval with coverage probability of at least $1 -
\alpha$.
We should also note once again that
the left endpoint of the interval may be larger than the right endpoint
in which case
the confidence interval
is taken to be the empty set, and so in this case
the length of this confidence interval is $\max( \delta_{j} - \delta
_{j+1} + ( \sqrt5 + \sqrt2 )z_{\alpha/2} \sigma_{j},0)$.

These results, for which a more formal proof is given in Section \ref
{proof.sec} are summarized in the following proposition.
%
\begin{proposition}
\label{m.CIj.prop}
For every $j \ge2$, the confidence interval $\mathrm{CI}_j^m$ defined in (\ref{mCIj})
has coverage probability of at least $1- {\alpha}$ for all
monotonically nondecreasing functions $f\in F_m$,
and for every $j \ge1$, the confidence interval $\mathrm{CI}_j^c$ defined in
(\ref{CIsj}) has
coverage probability of at least $1- {\alpha}$ for all convex functions
$f\in F_c$.
\end{proposition}

The second stage in the construction is that of
selecting from these collections of intervals the one to be used.
First note that one should not select the shortest interval since the
collections
defined in
\eqref{mCIj}
and
\eqref{CIsj}
will always contain one which corresponds to the empty set.
A more sensible goal is to try to select the interval with the smallest
expected length or at least one which has expected length close to
the smallest expected length.




The approach we take here
is to
choose an interval for which the expected length is of the same order
of magnitude
as the standard deviation of the length.
Such an interval will always have expected length close to the shortest
expected length.
For the case of monotonically nondecreasing functions the selection of
the interval from the countable collection
in (\ref{mCIj}) can be done by creating another collection of
estimators which
can be used to estimate the expected length of the intervals.

More specifically
set $\xi
_j=2^{j-1}(Y(2^{-j+1})-Y(2^{-j}))-2^{j-1}(Y(-2^{-j})-Y(-2^{-j+1}))$.
Then for $j\geq2$, $\xi_j$'s are independent of each other and both
$\delta_j^R$ and $\delta_j^L$ are independent of $\xi_k$ for every
$k\leq j$.
We should note that the estimators $\xi_j$ are similar to $\delta_j^R -
\delta_j^L$ in that they
are both differences of averages of $Y$ to the left and right of the
origin and thus estimate the average local change
of the function.
However $\delta_j^R - \delta_j^L$ are not independent for different $j$
whereas the $\xi_j$ are independent. It is thus natural to view the
$\xi
_j$ as a surrogate for $\delta_j^R- \delta_j^L$
with the technical advantage that they are independent.
The selection of a $j$ for which $\xi_j$ has expected value close to
$\sigma_j$ will then result in a confidence interval
$\mathrm{CI}_j^m$ close to the one with the smallest expected length.
The independence properties of the $\xi_j$ allows us to guarantee a $1-
\alpha$ coverage probability while making this selection.\vadjust{\goodbreak}

More specifically the construction proceeds as follows. Let
%
\begin{equation}
\label{mhat_j} \hat j = \inf_{j} \biggl\{ j \dvtx \xi_j \le\frac{3}{2}z_{\alpha} \sigma_j
\biggr\}
\end{equation}
and define the final confidence interval by
%
\begin{equation}
\label{m_CI} \mathrm{CI}^m_* = \mathrm{CI}_{\hat j}^m(
\alpha).
\end{equation}

Before we turn to the analysis of this procedure we also introduce here
a related confidence procedure in the convex case.
Here rather than introducing an independent estimate of the difference
between the two
estimators used in constructing the confidence interval, we proceed
more directly.
The basic idea is similar, but the dependence between the estimates of
$j$ and the
confidence interval constructed from this estimate requires that we
adjust the
original coverage level of our $\mathrm{CI}_j^c$.


More specifically
let $T_j = \delta_j - \delta_{j+1}$.
When the expected value of $T_j$ is the same order as $\sigma_j$, the
confidence interval $\mathrm{CI}_j^c$\vspace*{1pt}
will then be close to the one with the smallest expected length.
Our estimate of $j$ is given by an empirical version, namely
%
\begin{equation}
\label{hat.j} \hat j = \inf_{j} \{ j \dvtx T_j \le
z_{\alpha} \sigma_j\}.
\end{equation}

Although this estimate can be used to select the appropriate $\mathrm{CI}_j^c$
to use, as just mentioned,
care also needs to be taken to make sure that the resulting selected
interval maintains
the required coverage probability.
The analysis given below shows that a choice of $\alpha/6$ in the
construction of the original collection of intervals
guarantees an overall coverage probability of $\alpha$.
Thus in the case of convex functions, we define our interval by
%
\begin{equation}
\label{convex.CI} \mathrm{CI}^c_* = \mathrm{CI}_{\hat j}^c\biggl(
\frac{\alpha}{6}\biggr).
\end{equation}


\subsection{Analysis of the confidence intervals}
\label{analysis.sec}

In this section we present the properties of the confidence intervals
$\mathrm{CI}_*^m$ and $\mathrm{CI}_*^c$
defined by \eqref{m_CI} and \eqref{convex.CI} focusing on the coverage
and the expected length of these
intervals.

We begin with the confidence interval $\mathrm{CI}_*^m$.
In this case it is easy to check the coverage probability of $\mathrm{CI}_*^m$
by the independence of the interval $\mathrm{CI}_j^m$ and $\xi_k$ for every $k$
satisfying $2\leq k\leq j$.

The key to the analysis of the expected length is the introduction of
$j_*^m$ where
%
\begin{equation}
\label{mj*} j_*^m = \argmin_j \{j\dvtx E\xi_j
\le z_{\alpha}\sigma_j \}.
\end{equation}

The analysis of the expected length
relies on showing that $\hat j$ is highly concentrated around $j_*^m$.
The concentration of $\hat j$ around $j_*^m$ then provides
a bound on the expected length of $\mathrm{CI}^*$.
These results, for which a proof is given in Section~\ref{proof.sec}
are summarized in the following theorem.

\begin{theorem}
\label{m.CI-EL.thm}
Let $0<\alpha\leq0.2$. The confidence interval $\mathrm{CI}_*^m$ defined in
(\ref{m_CI})
has coverage probability of at least $1- \alpha$ for all
monotonically nondecreasing functions $f\in F_m$ and satisfies
%
\begin{eqnarray}
\label{b7} E_fL\bigl(\mathrm{CI}_*^m\bigr)\leq1.21(3z_{\alpha}+2
\sqrt{2}z_{{\alpha}/{2}})\sigma _{j_*^m}\leq c_0z_{\alpha}
\sigma_{j_*^m},
\end{eqnarray}
where $c_0$ is a constant and can be taken to be $8.85$ for all $0<
\alpha\le0.2$.
\end{theorem}

\begin{remark}\label{rem1}
The constant $c_0$ in Theorem~\ref{m.CI-EL.thm} depends on the upper
limit of $\alpha$. $c_0$ can be smaller if the upper limit on $\alpha$
is reduced. For example, for common choices of $\alpha=0.05$ or 0.01,
$c_0\le7.71$ for $\alpha= 0.05$, and $c_0\le7.42$ for $\alpha= 0.01$.
\end{remark}

Theorem~\ref{m.CI-EL.thm} shows that the coverage probability is
attained and also
provides an upper bound on the expected length in terms of $\sigma_{j_*^m}$.
In order to establish that this expected length is within a
constant factor of the lower bound given in
Theorem~\ref{EL.bound.thm}, we need to
provide a lower bound for $L_{\alpha}^{*}(f, F_m)$ in terms of $
z_{\alpha}\sigma_{j_*^m}$.
This connection is given in the following theorem.
%
\begin{theorem}
\label{mEL_lb}
Let $0 < \alpha\leq0.2$ and let $f \in F_m$. Then
%
\begin{equation}
\label{mEL_lb1} L_{\alpha}^{*}(f,
F_m) \ge\biggl(1 - \frac{1}{\sqrt{2 \pi}z_{\alpha}}\biggr) \frac
{1}{\sqrt{2}}z_{\alpha}
\sigma_{j_*^m}.
\end{equation}
\end{theorem}
Combining Theorems~\ref{m.CI-EL.thm} and~\ref{mEL_lb}, we have
%
\begin{equation}
E_f L\bigl(\mathrm{CI}_*^m\bigr) \le c_1
L_{\alpha}^{*}(f, F_m)
\end{equation}
for all monotonically nondecreasing functions $f\in F_m$, where $c_1$
is a constant depending on $\alpha$ only. For example, $c_1$ can be
taken to be $14.40$ for $\alpha=0.05$ and~$12.67$ for $\alpha=0.01$.
Hence, the confidence interval $\mathrm{CI}_*^m$ is uniformly within a constant
factor of the benchmark $L_{\alpha}^{*}(f, F_m) $ for all monotonically
nondecreasing functions $f$ and all confidence level $1-\alpha\ge0.8$.

We now turn to an analysis of the properties of the confidence interval
$\mathrm{CI}_*^c$ defined in \eqref{convex.CI}.
The key to this analysis is the introduction of $j_*^c$ where
%
\begin{equation}
\label{j*} j_*^c = \argmin_j \biggl\{j\dvtx ET_j
\le\frac{ 2}{3}z_{\alpha}\sigma_j \biggr\}.
\end{equation}

The analysis of both the coverage probability and the expected length
relies on showing that $\hat j$ is highly concentrated around $j_*^c$.
The probability of not covering $f(0)$ can be bounded
by
%
\begin{eqnarray}
P\bigl( f(0) \notin \mathrm{CI}_*^c\bigr) &\le& P\bigl(\hat j \le
j_*^c -3\bigr) + P\bigl(\hat j \ge j_*^c +3\bigr)
\nonumber
\\[-8pt]
\\[-8pt]
\nonumber
&&{} + \sum
_{l=-2}^{2} P\bigl( f(0) \notin
\mathrm{CI}_{j_*^c + l}\bigr).
\end{eqnarray}
The first two terms are controlled by the high concentration of $\hat
j$ around $j_*^c$,
and the last term is controlled by Proposition~\ref{m.CIj.prop} which
bounds the\vadjust{\goodbreak} coverage probability of any given $j$.
The concentration of $\hat j$ around $j_*^c$ also allows control on
the expected length of $\mathrm{CI}_*^c$ which leads to
the following theorem.


\begin{theorem}
\label{convex.CI-EL.thm}
Let $0<\alpha\leq0.2$. The confidence interval $\mathrm{CI}_*^c$ defined in
(\ref{convex.CI})
has coverage probability of at least $1- \alpha$ for all
convex $f$ and satisfies
%
\begin{equation}
E_f L\bigl(\mathrm{CI}_*^c\bigr) \le1.25 \bigl(z_\alpha+
(\sqrt{5} + \sqrt{2}) z_{\alpha /12}\bigr) \sigma_{j_*^c} \le
c_0 z_{\alpha}\sigma_{j_*^c},
\end{equation}
where $c_0$ is a constant and can be taken to be $12.79$ for all $0<
\alpha\le0.2$.
\end{theorem}

\begin{remark}\label{rem2}
The constant $c_0$ in Theorem~\ref{convex.CI-EL.thm} depends on the
upper limit of $\alpha$. $c_0$ can be smaller if the upper limit on
$\alpha$ is reduced. For example, for common choices of $\alpha=0.05$
or 0.01, $c_0\le8.57$ for $\alpha= 0.05$, and $c_0\le7.42$ for
$\alpha= 0.01$.
\end{remark}

Theorem~\ref{convex.CI-EL.thm} shows that the coverage probability is
attained and also
provides an upper bound on the expected length in terms of $\sigma_{j_*^c}$.
As was the case for monotone functions,
in order to to establish that this expected length for convex functions
is within a
constant factor of the lower bound given in
Theorem~\ref{EL.bound.thm}, we need to
provide a lower bound for $L_{\alpha}^{*}(f, F_c)$ in terms of $
z_{\alpha}\sigma_{j_*^c}$.
This connection is given in the following theorem.
%
\begin{theorem}
\label{EL.bound2.thm}
Let $0 < \alpha\leq0.2$ and let $f \in F_c$. Then
%
\begin{equation}
\label{EL.bound3} L_{\alpha}^{*}(f, F_c) \ge
\biggl(1 - \frac{1}{\sqrt{2 \pi}z_{\alpha}}\biggr) \frac
{\sqrt2}{3}z_{\alpha}
\sigma_{j_*^c}.
\end{equation}
\end{theorem}
Theorems~\ref{convex.CI-EL.thm} and~\ref{EL.bound2.thm} together yield
%
\begin{equation}
E_f L\bigl(\mathrm{CI}_*^c\bigr) \le c_2
L_{\alpha}^{*}(f, F_c)
\end{equation}
for all convex functions $f\in F_c$, where $c_2$ is a constant
depending on $\alpha$ only. For example, $c_2$ can be taken to be $24$
for $\alpha=0.05$ and $19$ for $\alpha=0.01$. Hence, the confidence
interval $\mathrm{CI}_*^c$ is uniformly within a constant factor of the
benchmark $L_{\alpha}^{*}(f, F_c) $ for all convex functions $f$ and
all confidence levels $1-\alpha\ge0.8$.




\section{Nonparametric regression}
\label{regression.sec}

We have so far focused on the white noise model. The theory presented
in the earlier sections can also easily be extended to nonparametric regression.
Consider the regression model
%
\begin{equation}
y_i = f(x_i) + \sigma z_i,\qquad i = -n,
-(n-1), -1, 0, 1,\ldots, n, \label{reg.model}
\end{equation}
where $x_i = {i\over2n}$ and $z_i \stackrel{\mathrm{i.i.d.}}{\sim} N(0, 1)$ and
where for notational convenience we index the observations from $-n$ to
$n$. Note that the noise level $\sigma$ can be accurately estimated
easily, as in Hall, Kay and
Titterington (\citeyear{HalKayTit90}) or \citet{MBWF}. See also \citet{Wanetal08}. We shall thus assume it is known in this section. Then
under the assumption that $f$ is convex or monotone, we wish to provide
a confidence interval for $f(0)$.


\subsection{Monotone regression}
Let $J=\lfloor\log_2 n\rfloor$. For $1 \le j \le J$ define the local
average estimators
%
\begin{equation}
\bar\delta_j^R = 2^{-j+1} \sum
_{k=1}^{2^{j-1}} y_{k} \quad\mbox{and}\quad \bar
\delta_j^L = 2^{-j+1} \sum
_{k=1}^{2^{j-1}} y_{ -k}.
\end{equation}
We should note that the indexing scheme is the reverse of that given
for the white noise with drift process.
Here estimators $\bar\delta_j^R$ (or $\bar\delta_j^L$) with small
values of $j$ have smaller bias (or larger bias) and larger variance
than those with larger values of $j$.

As in the white noise model it is easy to check that $\bar\delta_j^R$
has nonnegative bias and $\bar\delta_j^L$ has nonpositive bias. Simple
calculations show that the variance of $\bar{\delta}_j^R$ and $\bar
{\delta}_j^L$ are both $2\sigma_j^2$, where $\sigma_j^2=2^{-j}\sigma^2$.
It is also\vspace*{-1pt} important to introduce
$\bar{\xi}_j$
as in the white noise case, where $\bar{\xi}_j=2^{-j}\sum_{k=2^{j-1}+1}^{2^{j}}(y_k-y_{-k})
$. 
It is easy to check that $E\bar{\xi}_j\leq E\bar{\xi}_{j+1}$, $\bar
{\xi
}_j$'s are independent with each other, and both $\bar{\delta}_j^R$ and
$\bar{\delta}_j^L$ are independent with $\bar{\xi}_k$ for every
$k\geq j$.

It then follows that $\mathrm{CI}_j^m=[\bar{\delta}_j^L-z_{{\alpha}/{2}}\sqrt {2}\sigma_j,\bar{\delta}_j^R+z_{{\alpha}/{2}}\sqrt{2}\sigma
_j]$ has
guaranteed coverage probability of at least $1- \alpha$ over all
monotonically nondecreasing functions.

Now
set
\begin{equation}
\label{reg.hat.j} \hat j = \cases{ %
\displaystyle\max
_{j} \biggl\{ j \dvtx \bar{\xi}_j \le
\frac{3}{2}z_{\alpha}\sigma _j \biggr\}, &\quad $\mbox{if $\displaystyle
\bar{\xi}_1 \le\frac{3}{2}z_{\alpha
}\sigma
_1$;}$
\vspace*{2pt}\cr
1, &\quad $\mbox{otherwise,}$}
\end{equation}
and define the confidence interval to
be
%
\begin{equation}
\label{reg.mCI} \mathrm{CI}_*^m = \mathrm{CI}_{\hat j}^m.
\end{equation}
The properties of this confidence interval can then be analyzed in the
same way as before and can be shown to be similar to those for the
white noise model. In particular, the following result holds.
%
\begin{theorem}
\label{mEL_reg.thm}
Let $0<\alpha\le0.2$. The confidence interval $\mathrm{CI}_*^m $ defined in
\eqref{reg.mCI} has coverage probability of at least $1- \alpha$ for all
monotone functions $f$ and satisfies
%
\begin{equation}
E_f L\bigl(\mathrm{CI}_*^m\bigr) \le C_1
L_{\alpha}^{*}(f, F_m)
\end{equation}
for all monotonically nondecreasing functions $f\in F_m$, where $C_1>0$
is a constant depending on $\alpha$ only.
\end{theorem}


\subsection{Convex regression}
As in the monotone case, set $J=\lfloor\log_2 n\rfloor$. For $1 \le j
\le J$ define the local average estimators
%
\begin{equation}
\bar\delta_j = 2^{-j} \sum_{k=1}^{2^{j-1}}
(y_{ -k} + y_{k}).\vadjust{\goodbreak}
\end{equation}
%
We should note that this indexing scheme is the reverse of that given
for the white noise with drift process.
Here estimators $\bar\delta_j$ with small values of $j$ have smaller
bias and larger variance
than those with larger values of $j$. 

As in the white noise model it is easy to check that $\bar\delta_j$
has nonnegative bias.
It is also important to introduce an estimate which has a similar
variance but is guaranteed to
have nonpositive bias. The key step is to introduce
%
\begin{equation}
\label{reg.Tj} T_j = \bar\delta_{j} - \bar
\delta_{j-1}
\end{equation}
as an estimate of the bias of $\bar\delta_j$.
The following lemma gives the required properties of $\bar\delta_j$
and $T_j$.
%
\begin{lemma}
\label{reg.bias.lem}
For any convex function $f$,
%
\begin{eqnarray}
2 E T_{j} &\le& E T_{j+1},\label{reg.bias.bound1}
\\
0 &\le&\operatorname{ Bias}(\bar\delta_{j})\le {2^{j-1}+1\over2^{j}+1} \operatorname{
Bias}(\bar\delta_{j+1}). \label{reg.bias.bound2} 
\end{eqnarray}
%
\end{lemma}
From \eqref{reg.bias.bound2} it is clear that the biases of
the estimators $\bar\delta_{j}$ are nonnegative and monotonically
nondecreasing. In addition straightforward
calculations using both~\eqref{reg.bias.bound1} and \eqref
{reg.bias.bound2} show that the estimators
\[
\delta^L_j = \bigl(2+2^{-(j-1)}\bigr)\bar
\delta_{j} - \bigl(1+2^{-(j-1)}\bigr)\bar \delta _{j+1}
= \bar\delta_{j} - \bigl(1+2^{-(j-1)}\bigr) T_{j+1}
\]
have a nonpositive and monotonically nonincreasing biases.
Simple calculations show that the variance of $\delta^L_j$
is $\tau_j^2 = (5+2^{-j+3} + 2^{-2j+2}) 2^{-j-1} \sigma^2$.

It then follows that
$\mathrm{CI}_j^c = [\bar\delta_j - (1+2^{-(j-1)}) T_{j+1} - z_{\alpha/12}
\tau
_j, \bar\delta_j + z_{\alpha/12} \sigma_j]$
has coverage over all convex functions.

Now
set
%
\begin{equation}
\label{c.reg.hat.j} \hat j = \cases{ %
\displaystyle\max
_{j} \{ j \dvtx T_j \le z_{\alpha}
\sigma_j \}, & \quad$\mbox{if $T_2 \le z_{\alpha}
\sigma_2$;}$
\vspace*{2pt}\cr
1, & \quad$\mbox{otherwise,}$ }
\end{equation}
and define the confidence interval to
be
%
\begin{equation}
\label{reg.cCI} \mathrm{CI}_*^c = \mathrm{CI}_{\hat j}^c.
\end{equation}
This confidence interval shares similar properties as the one for the
white noise model. In particular, the following result holds.
%
\begin{theorem}
\label{cEL_reg.thm}
Let $0<\alpha\le0.2$. The confidence interval $\mathrm{CI}_*^c$ defined in
\eqref{reg.cCI} has coverage probability of at least $1- \alpha$ for all
convex function $f$ and satisfies
%
\begin{equation}
E_f L\bigl(\mathrm{CI}_*^c\bigr) \le C_2
L_{\alpha}^{*}(f, F_c)
\end{equation}
for all convex function $f\in F_c$, where $C_2>0$ is a constant
depending on $\alpha$ only.
\end{theorem}

\section{Discussion}
\label{discussion.sec}

The major emphasis of the paper has been to show that with shape
constraints it is
possible to construct confidence intervals that have expected length
that adapts to individual functions.
In this section we shall discuss briefly
the maximum expected lengths of our procedures over Lipschitz classes
that are
either monotone or convex in a way that
is similar to that provided in
D\"{u}mbgen (\citeyear{Dumbgen,Dum03}) for the maximum width of a
confidence band.
We
shall also explain how our results can be extended to the problem of
estimating the value of $f$ at points other than $0$.

\subsection{Minimax results}

Although the focus of the present paper has been on the construction of
a confidence interval
with the expected length adaptive to each individual convex or monotone
function, these results
do yield immediately adaptive minimax results for the expected length
in the conventional sense.
Define
\[
F_c(\beta, M) = F_c\cap\Lambda(\beta, M)\quad \mbox{and}\quad
F_m(\beta, M) = F_m\cap\Lambda(\beta, M).
\]
The following results are direct consequence of Theorems \ref
{m.CI-EL.thm} and~\ref{convex.CI-EL.thm}.
%
\begin{corollary}
\textup{(i)} The confidence interval $\mathrm{CI}_*^m$ defined in (\ref{m_CI})
satisfies
%
\begin{equation}
\sup_{f \in F_m(\beta, M)} E_fL\bigl(\mathrm{CI}_*^m\bigr)
\le C_1 M^{1/(1+2\beta)}n^{-{\beta/(1 + 2 \beta)}}
\end{equation}
simultaneously for all $0 \le\beta\le1$ and $1 < M < \infty$, for
some absolute constant $C_1>0$.

\textup{(ii)} The confidence interval $\mathrm{CI}_*^c$ defined in (\ref{convex.CI})
satisfies
%
\begin{equation}
\sup_{f \in F_c(\beta, M)} E_fL\bigl(\mathrm{CI}_*^c\bigr)
\le C_2 M^{1/(1+2\beta)}n^{-{\beta/(1 + 2 \beta)}}
\end{equation}
simultaneously for all $1\le\beta\le2$ and $1 < M < \infty$, for
some absolute constant $C_2>0$.
\end{corollary}

We should note that these ranges of Lipschitz classes are the only ones
of interest in these cases.
In particular
suppose that $\mathrm{CI}$ is a confidence interval with guaranteed coverage
over the class
of monotonically nondecreasing functions.
Then for
any $\beta> 1$ the class $\Lambda(\beta,M)$ includes the linear
function $f_k(t) = kt$.
As shown in Example~\ref{linear.ex} in Section~\ref{example.sec},
\[
L^*_\alpha(f_k, F_m) \ge\biggl( 1 -
\frac{1}{\sqrt{2\pi} z_{\alpha}}\biggr) (3k)^{{1}/{3}}z^{2/3}_{\alpha}n^{-{1/3}}.
\]
Hence,
%
\begin{eqnarray}
\sup_{f \in F_m(\beta, M)} E_fL(\mathrm{CI}) &\ge&\sup_{k}L^*_\alpha(f_k,
F_m) \nonumber
\\[-8pt]
\\[-8pt]
\nonumber
&=& \sup_{k}\biggl( 1 - \frac{1}{\sqrt{2\pi} z_{\alpha}}
\biggr) (3k)^{{1}/{3}}z^{2/3}_{\alpha}n^{-{1/3}}
=
\infty.
\end{eqnarray}

A similar results holds for convex functions assumed to belong to
$\Lambda(\beta,M)$ with $\beta>2$.
On the other hand suppose $f$ is convex and assumed to belong to
$\Lambda(\beta,M)$ with $\beta<1$.
Then from the assumption that $f$ is in $\Lambda(\beta,M)$
it follows that $|f(1/2) -f(-1/2)| \le M$. Convexity then shows that $f
\in\Lambda(1,M)$
and the maximum expected length over this class is given above.

\subsection{Confidence interval at other points}

The focus of the present paper has been on the problem
of estimating the value of $f(0)$.
The basic development is similar for any other point $t$ in the
interior of the interval $[-1/2,1/2]$
unless $t$ is near to the boundary.
More specifically for any $0 \le t < 1/2$
we can consider estimators $\delta_j^R(t)=2^j(Y(t+2^{-j})-Y(t))$ and
$\delta_j^L(t)=2^j(Y(t)-Y(t-2^{-j}))$ where $j \ge-\log_2(\frac
{1}{4}-\frac{t}{2})$ for monotone functions and
$\delta_j(t) = 2^{j-1}(Y(t +2^{-j}) - Y(t-2^{-j}))$
where $j \ge-\log_2(\frac{1}{2}-t)$ for convex functions.
The basic theory is the same as before.

For monotonically nondecreasing functions, the confidence interval
$\mathrm{CI}_j^m$ is replaced by
\[
\mathrm{CI}_j^m(t) = \bigl[\delta_j^L(t)-z_{{\alpha}/{2}}
\sqrt{2}\sigma _j,\delta _j^R(t)+z_{{\alpha}/{2}}
\sqrt{2}\sigma_j\bigr]
\]
and the choice of $\hat j$ is given by
\[
\hat j(t) = \inf_{j \ge-\log_2({1}/{4}-({t}/{2}))} \biggl\{ j \dvtx \xi_j(t) \le
\frac{3}{2}z_{\alpha} \sigma_j \biggr\},
\]
where $\xi
_j(t)=2^{j-1}(Y(t+2^{-j+1})-Y(t+2^{-j}))-2^{j-1}(Y(t-2^{-j})-Y(t-2^{-j+1}))$.
The final confidence interval is defined by
%
\begin{equation}
\label{m.CI.t} \mathrm{CI}_*^m = \mathrm{CI}^m_{\hat j(t)}.
\end{equation}

For convex functions, the confidence interval $\mathrm{CI}_j^c$ is replaced by
\[
\mathrm{CI}_j^c(t) = \bigl[ \delta_{j+1}(t) - \bigl(
\delta_j(t) - \delta_{j+1}(t)\bigr)_+ - z_{\alpha/12}
\sqrt{5}\sigma_j, \delta_{j+1}(t) + z_{\alpha/12}
\sigma_{j+1}\bigr],
\]
and $\hat j$ is chosen to be
\[
\hat j(t) = \inf_{j \ge-\log_2({1}/{2}-t)} \bigl\{ j \dvtx T_j(t) \le
z_{\alpha} \sigma_j \bigr\}.
\]
Define the final confidence interval by
\[
\mathrm{CI}_*^c = \mathrm{CI}^c_{\hat j(t)}.
\]

The modulus of continuity defined in \eqref{l.modulus} is replaced by
%
\begin{equation}
\label{l.modulus} \omega(\ep, f, t,\mathcal{F}) = \sup\bigl\{ \bigl|g(t) - f(t)\bigr| \dvtx \|g - f \|_2 \le\ep, g \in\mathcal{F}\bigr\}.
\end{equation}
The earlier analysis then yields
\[
E_f L\bigl(\mathrm{CI}_*^m(t)\bigr) \le c_1
L^*_{\alpha}(f,t,F_m)
\]
and
\[
E_f L\bigl(\mathrm{CI}_*^c(t)\bigr) \le c_2
L^*_{\alpha}(f,t,F_c),
\]
where we now have
\[
L^*_\alpha(f,t,\mathcal{F}) \ge\biggl( 1 - \frac{1}{\sqrt{2\pi}
z_{\alpha}}\biggr)
\omega\biggl(\frac{z_{\alpha}}{\sqrt n}, f,t,\mathcal{F}\biggr).
\]


Finally we should note that at the boundary the
construction of a confidence interval must be unbounded. For example
any honest confidence interval for $f(1/2)$ must be of the form $[\hat
f(1/2), \infty)$; otherwise it cannot have
guaranteed coverage probability.


%
%







\section{Proofs}
\label{proof.sec}
We prove the main results in this section. We shall omit the proofs for
Theorems~\ref{mEL_reg.thm} and~\ref{cEL_reg.thm} as they are analogous
to those for the corresponding results in the white noise model.


\subsection{\texorpdfstring{Proof of Lemma \protect\ref{con.lem}}{Proof of Lemma 1}}
Set $f_s(t) = \frac{f(t) + f(-t)}{2} - f(0)$.
Now note that $f_s(tx)$ is convex in $x$ for all $0 \le t \le1$.
Hence $g(x) = \int_{0}^{1} f_s(tx) \,dt$ is also convex with $g(0) = 0$.
For $x > 0$ set $z =xt$, and it follows that
$
g(x) = \frac{1}{x} \int_{0}^{x} f_s(z) \,dz =
\frac{1}{2x} \int_{-x}^{x} (f(z) - f(0)) \,dz$.
Equation (\ref{con.bias.eq}) follows from the fact that
$g(x) \le\frac{1}{2} g(2x)$ for $x = 2^{-(j+1)}$, and equation (\ref
{con.eq}) follows from the fact that $g(2x)\leq2/3g(x)+1/3g(4x)$. 
\subsection{\texorpdfstring{Proof of Lemma \protect\ref{reg.bias.lem}}{Proof of Lemma 2}}
For any convex function $f$, let $f_s(x)=\frac{1}{2}(f(x)+f(-x)) - f(0)$.
Then $f_s(x)$ is convex, increasing in $|x|$ and $f_s(0)=0$. Convexity
of $f_s$ yields
that for $0<x\le y$,
%
\begin{equation}
\label{cnvx} {f_s(x)\over x} \le{f_s(y)\over y}.
\end{equation}
Note that $E \delta_j = 2^{-(j-1)}\sum_{k=1}^{2^{j-1}} f_s({k\over
n})$ and
\[
E T_j = 2^{-(j-1)} \Biggl\{\sum_{k=2^{j-2}+1}^{2^{j-1}}
f_s\biggl({k\over n}\biggr) - \sum
_{k=1}^{2^{j-2}} f_s\biggl(
{k\over n}\biggr) \Biggr\}.
\]
So $ET_j \ge2 ET_{j-1}$ is equivalent to
\[
\sum_{k=2^{j-2}+1}^{2^{j-1}} f_s\biggl(
{k\over n}\biggr) - \sum_{k=1}^{2^{j-2}}
f_s\biggl({k\over n}\biggr) \ge4 \sum
_{k=2^{j-3}+1}^{2^{j-2}} f_s\biggl(
{k\over n}\biggr) - 4\sum_{k=1}^{2^{j-3}}
f_s\biggl({k\over n}\biggr),
\]
which is the same as
%
\begin{equation}
\label{aaaaa} \sum_{k=2^{j-2}+1}^{2^{j-1}}
f_s\biggl({k\over n}\biggr) + 3\sum
_{k=1}^{2^{j-3}} f_s\biggl(
{k\over n}\biggr) \ge5 \sum_{k=2^{j-3}+1}^{2^{j-2}}
f_s\biggl({k\over n}\biggr).
\end{equation}
Now note that for $x \ge0$ and $u \ge0$,
\begin{eqnarray*}
f_s(x) + f_s(x + 3u) &\ge& f_s(x+u) +
f_s(x + 2u) \quad\mbox{and}\\
 f_s(x) + f_s(x+2u)
&\ge&2f_s(x+ u)
\end{eqnarray*}
and consequently
$
f_s(x+3u) + f_s(x+2u) + 3 f_s(x) \ge5f_s(x+u).
$
Then (\ref{aaaaa}) follows by taking $u={2^{j-3}\over n}$ and
$x={k\over n}$ and then summing over $k=1,\ldots, 2^{j-3}$.

Denote the bias of $\bar\delta_j$ by $\bar b_j = E\bar\delta_j -
f(0)$. Then
\[
\bar b_j = 2^{-(j-1)}\sum_{k=1}^{2^{j-1}}
f_s\biggl({k\over
n}\biggr)=2^{-(j-1)} \Biggl\{
\sum_{k=2^{j-2}+1}^{2^{j-1}} f_s\biggl(
{k\over n}\biggr) + \sum_{k=1}^{2^{j-2}}
f_s\biggl({k\over n}\biggr) \Biggr\}.
\]
It follows from (\ref{cnvx}) that for $k > 2^{j-2}$,
$f_s({k\over n}) \ge{k\over2^{j-2}} f_s({2^{j-2}\over n})$,
and for $k \le2^{j-2}$,
$f_s({k\over n}) \le{k\over2^{j-2}} f_s({2^{j-2}\over n})$. Hence
\begin{eqnarray*}
\sum_{k=2^{j-2}+1}^{2^{j-1}} f_s\biggl(
{k\over n}\biggr) &\ge& \sum_{k=2^{j-2}+1}^{2^{j-1}}
{k\over2^{j-2}} \cdot f_s\biggl({2^{j-2}\over n}
\biggr)\ge {\sum_{k=2^{j-2}+1}^{2^{j-1}}{k/2^{j-2}}\over\sum_{k=1}^{2^{j-2}}{k/2^{j-2}}} \sum_{k=1}^{2^{j-2}}f_s
\biggl({k\over n}\biggr)
\\
&=&{3\cdot2^{j-2} + 1\over2^{j-2}+1}\sum_{k=1}^{2^{j-2}}f_s
\biggl({k\over n}\biggr).
\end{eqnarray*}
Hence,
\[
\bar b_j \ge2^{-(j-1)}\cdot\biggl({3\cdot2^{j-2} + 1\over
2^{j-2}+1}+1
\biggr)\sum_{k=1}^{2^{j-2}}f_s\biggl(
{k\over n}\biggr) 
={2^{j-1} + 1\over2^{j-2}+1}
\bar b_{j-1}. 
\]

\subsection{\texorpdfstring{Proof of Theorem \protect\ref{EL.bound.thm}}{Proof of Theorem 1}}

Suppose that $X \sim N(\theta, \sigma^2)$ where it is known that
$\theta\in[0,a \sigma]$.
The confidence interval for $\theta$ which has guaranteed coverage
over the interval $\theta\in[0,a \sigma]$
and which minimizes the expected length when $\theta= 0$ is given by
%
\begin{equation}
\bigl[0, \max\bigl(0, \min(X + z_{\alpha} \sigma, a \sigma)\bigr)\bigr].
\end{equation}

It follows that
%
\begin{equation}
\qquad L = \sigma\int_{-z_{\alpha}}^{a - z_{\alpha}} z \phi(z) \,dz + \sigma
\bigl(z_{\alpha} P(-z_{\alpha} \le Z \le a - z_{\alpha}) + a
P(Z \ge a - z_{\alpha})\bigr)
\end{equation}
and hence
%
\begin{eqnarray}
\frac{L}{\sigma} &=& \bigl(\phi(z_{\alpha}) - \phi(a - z_{\alpha})
\bigr) + z_{\alpha} \bigl(\Phi( a- z_{\alpha}) -
\Phi(-z_{\alpha})\bigr)
\nonumber
\\[-8pt]
\\[-8pt]
\nonumber
&&{}+ a \bigl(1- \Phi ( a - z_{\alpha})\bigr).
\end{eqnarray}


In particular when $a = z_{\alpha}$,
%
\begin{equation}
\frac{L}{\sigma} \ge z_{\alpha}\biggl( 1 - \frac{\phi(0)}{z_{\alpha}} +
\frac{\phi(z_{\alpha
})}{z_{\alpha}} - \alpha\biggr).
\end{equation}
In particular we have
%
\begin{equation}
\frac{L}{\sigma} \ge z_{\alpha}\biggl( 1 - \frac{\phi(0)}{z_{\alpha}} \biggr).
\end{equation}

Write $L_{\alpha}^*(f, \cal{F})$ for the smallest expected length at
$f$ when
we have guaranteed coverage over $\cal F$.
In particular let $P_{\theta}$ be a subfamily of $\cal F$, and then
$L_{\alpha}^*(f, {\cal F}) \ge L_{\alpha}^*(f, P_{\theta})$.

Now suppose that $f_0$ is the ``true'' function. 
Fix $\ep> 0$. There is a function $f_1 \in\sf$ such that
\[
\bigl| f_1(0) - f_0(0)\bigr | = \omega\biggl(
{\varepsilon\over{\sqrt n}}, f,\sf\biggr)
\]
and such that
\[
\| f_1 - f_0 \|_2 = {\varepsilon\over\sqrt n}.
\]

Now for $0 \le\theta\le1$, let $f_{\theta} = f_0 + \theta(f_1 - f_0)$.
Let $P_{\theta}$ be this collection of $f_{\theta}$.
Now for the process
\[
dY(t) = f_{\theta}(t) \,dt + {1\over\sqrt n} \,dW(t), \qquad -
{1\over2}\leq t\leq{1\over2},
\]
there is a sufficient statistic $S_n$ given by
\[
S_n = f_0(0) + \bigl(f_1(0) -
f_0(0)\bigr) \frac{1}{\int(f_1 - f_0)^2}\int \bigl(f_1(t)-
f_0(t)\bigr) \bigl(dY(t) - f_0(t) \,dt\bigr).
\]
Note that $S_n$ has
a normal distribution $S_n \sim N(f_{\theta}(0), \frac{(f_1(0) -
f_0(0))^2 }{n \int(f_1 - f_0)^2})$
or more specifically $S_n \sim N(f_{\theta}(0), \frac{1}{\varepsilon
^2}\omega^2(\frac{\varepsilon}{\sqrt{n}}, f_0,\sf))$.

Note that $a = \varepsilon$. Now take $\varepsilon= z_{\alpha}$. It then
follows that
\[
L_{\alpha}^*(f_0, P_{\theta}) \ge\omega\biggl(
\frac{z_{\alpha}}{\sqrt
n},f_0,\sf\biggr) \biggl( 1 - \frac{\phi(0)}{z_{\alpha}} +
\frac{\phi(z_{\alpha
})}{z_{\alpha}} - \alpha\biggr). 
\]

\subsection{\texorpdfstring{Proof of Proposition \protect\ref{m.CIj.prop}}{Proof of Proposition 1}}
For monotone functions, we have
\begin{eqnarray*}
\pr\bigl(f(0)\in \mathrm{CI}_j^m\bigr)&=& \pr\bigl(
\delta_j^L-z_{{\alpha}/{2}}\sqrt{2}\sigma_j
\leq f(0)\leq \delta _j^R+z_{{\alpha}/{2}}\sqrt{2}
\sigma_j\bigr)
\cr
&\geq&1-\pr\bigl(\delta_j^R<f(0)-z_{{\alpha}/{2}}
\sqrt{2}\sigma _j\bigr)-\pr \bigl(\delta_j^L>f(0)+z_{{\alpha}/{2}}
\sqrt{2}\sigma_j\bigr)
\cr
&=&1-\pr\biggl(Z<\frac{f(0)-E(\delta_j^R)}{\sqrt{2}\sigma_j}-z_{{\alpha
}/{2}}
\biggr)\\
&&{}-\pr\biggl(Z>\frac{f(0)-E(\delta_j^L)}{\sqrt{2}\sigma_j}+z_{
{\alpha}/{2}}\biggr),
\end{eqnarray*}
where $Z$ is a standard normal random variable. Because $f(0)-E(\delta
_j^R)\leq0$ and $f(0)-E(\delta_j^L)\geq0$, we have
\[
\pr\bigl(f(0)\in \mathrm{CI}_j^m\bigr)\geq 1-
\pr(Z<-z_{{\alpha}/{2}})-\pr (Z>z_{{\alpha}/{2}})=1-\alpha.
\]
For convex functions, let $b_j=\operatorname{ Bias}(\delta_{j})$. It follows from
Lemma~\ref{con.lem} that $b_j-2b_{j+1}>0$, and hence
we have
\begin{eqnarray*}
P\bigl(f(0) \in \mathrm{CI}_j^c\bigr) &\ge& 
P\bigl(2
\delta_{j+1}-\delta_j - z_{{\alpha}/{2}} \sqrt{5}\sigma
_{j}\le f(0)\le\delta_{j+1} + z_{{\alpha}/{2}}
\sigma_{j+1}\bigr)
\\
&\ge& 1- P\bigl( \delta_{j+1} <f(0) -
z_{{\alpha}/{2}}\sigma_{j+1}\bigr) \\
&&{}- P\bigl(2\delta_{j+1}-
\delta_j > f(0) + z_{{\alpha}/{2}}\sqrt {5}\sigma _{j}
\bigr)
\\
&=& 1- P\biggl( {\delta_{j+1} -E\delta_{j+1}\over\sigma
_{j+1}}<-{b_{j+1}\over\sigma_{j+1}}-
z_{{\alpha}/{2}}\biggr)
\\
&&{}- P\biggl({2\delta_{j+1}-\delta_j - E(2\delta_{j+1}-\delta_j )\over
\sqrt {5}\sigma_j} > {b_j - 2b_{j+1}\over\sqrt{5}\sigma_j}+
z_{{\alpha
}/{2}}\biggr)
\\
&\ge& 1 - P(Z<-z_{{\alpha}/{2}}) - P(Z > z_{{\alpha
}/{2}})
\\
& = & 1 - \alpha. 
\end{eqnarray*}

\subsection{\texorpdfstring{Proof of Theorem \protect\ref{m.CI-EL.thm}}{Proof of Theorem 2}}
We shall first prove that the confidence interval $\mathrm{CI}_*^m$ has
guaranteed coverage probability of $1-\alpha$ over $F_m$ and then prove
the upper bound for the expected length.

Note that
\[
\pr\bigl(f(0)\in \mathrm{CI}_*^m\bigr)=\sum_{j=2}^{\infty}
\pr\bigl(f(0)\in \mathrm{CI}_j^m | \hat {j}=j\bigr)\pr(\hat{j}=j).
\]
Because both $\delta_j^R$ and $\delta_j^L$ are independent of $\xi_k$
for $k\leq j$, and the event $\{\hat{j}=j\}$ depends only on $\xi_k$
for $k\leq j$, then by Proposition~\ref{m.CIj.prop} we have
\[
\pr\bigl(f(0)\in \mathrm{CI}_*^m\bigr)=\sum_{j=2}^{\infty}
\pr\bigl(f(0)\in \mathrm{CI}_j^m \bigr)\pr (\hat {j}=j)\geq\sum
_{j=2}^{\infty}(1-\alpha)\pr(\hat{j}=j)=1-
\alpha.
\]

We now turn to the upper bound for the expected length.
Note that for $s\geq0$, $E\xi_{j_*^m+s}\leq z_{\alpha}\sigma
_{j_*^m}=\frac{1}{2^{{s}/{2}}}z_{\alpha}\sigma_{j_*^m+s}$, and so
we have
\begin{eqnarray*}
\pr\bigl(\hat{j}\geq j_*^m+k\bigr)&\leq&\prod
_{s=0}^{k-1}\pr\biggl(\xi _{j_*^m+s}>
\frac
{3}{2}z_{\alpha}\sigma_{j_*^m+s}\biggr)\\
&\leq&\prod
_{s=0}^{k-1}\pr \biggl(Z>z_{\alpha
}\biggl(
\frac{3}{2}-\frac{1}{2^{{s}/{2}}}\biggr)\biggr).
\end{eqnarray*}
It follows from $E(\delta_j^R-\delta_j^L)\leq2E\xi_j$ that
$E(\delta
_{\hat{j}}^R-\delta_{\hat{j}}^L)\leq2E\xi_{\hat{j}}$, and hence
we have
%
\begin{eqnarray*}
E_fL\bigl(\mathrm{CI}^*\bigr)&=&E_f\bigl(\delta_{\hat{j}}^R-
\delta_{\hat{j}}^L+2\sqrt {2}z_{
{\alpha}/{2}}
\sigma_{\hat{j}}\bigr)\leq E_f(2\xi_{\hat{j}}+2\sqrt
{2}z_{
{\alpha}/{2}}\sigma_{\hat{j}})
\cr
&\leq&E_f
\bigl((3z_{\alpha}+2\sqrt{2}z_{{\alpha}/{2}})\sigma_{\hat
{j}}\bigr)=
\sum_{j=2}^{\infty}(3z_{\alpha}+2
\sqrt{2}z_{{\alpha
}/{2}})\sigma_{j}\cdot\pr(\hat{j}=j).
\end{eqnarray*}
Thus
%
\begin{equation}
\label{b8}  E_fL\bigl(\mathrm{CI}^*\bigr)\leq(3z_{\alpha}+2
\sqrt{2}z_{{\alpha}/{2}})\sigma _{j_*^m} {\Biggl(}\pr\bigl(\hat{j}\leq
j_*^m\bigr)+\sum_{k=1}^{\infty
}2^{k/2}
\pr\bigl(\hat {j}=j_*^m+k\bigr) {\Biggr)}.\hspace*{-35pt}
\end{equation}
Set $w_k=2^{k/2}-2^{(k-1)/2}$ for $k\geq1$. Then it is easy to see that
\[
S=\pr\bigl(\hat{j}\leq j_*^m\bigr)+\sum
_{k=1}^{\infty}2^{k/2}\pr\bigl(\hat
{j}=j_*^m+k\bigr)=1+\sum_{k=1}^{\infty}w_k
\pr\bigl(\hat{j}\geq j_*^m+k\bigr).
\]
Thus
\[
S=1+\sum_{k=1}^{\infty}w_k\prod
_{s=0}^{k-1}\pr\biggl(Z>z_{\alpha}
\biggl(\frac
{3}{2}-\frac{1}{2^{{s}/{2}}}\biggr)\biggr).
\]
The right-hand side is increasing in $\alpha$. Through numerical
calculations, we can see that, for $\alpha=0.2$,
\[
\sum_{k=1}^{\infty}w_k\prod
_{s=0}^{k-1}\pr\biggl(Z>z_{\alpha}
\biggl(\frac
{3}{2}-\frac{1}{2^{{s}/{2}}}\biggr)\biggr)\leq0.21.
\]
Thus, by equation (\ref{b8}), we have
\[
E_fL\bigl(\mathrm{CI}^*\bigr)\leq1.21(3z_{\alpha}+2
\sqrt{2}z_{{\alpha
}/{2}})\sigma _{j_*^m}. 
\]

\subsection{\texorpdfstring{Proof of Theorem \protect\ref{mEL_lb}}{Proof of Theorem 3}}
Note that if $j_*^m>2$, then $E\xi_{j_*^m-1}\geq z_{\alpha}\sigma
_{j_*^m-1}=\frac{1}{\sqrt{2}}z_{\alpha}\sigma_{j_*^m}$ and hence there
is a $t_*\leq2^{-j_*^m+2}$ such that we have either $f(t_*)-f(0)\geq
\frac{1}{\sqrt{2}}z_{\alpha}\sigma_{j_*^m}$ or $f(0)-f(-t_*)\geq
\frac
{1}{\sqrt{2}}z_{\alpha}\sigma_{j_*^m}$. If $f(t_*)\geq\frac
{1}{\sqrt {2}}z_{\alpha}\sigma_{j_*^m}+f(0)$
, let
\[
g(t)=\cases{ %
\displaystyle\max\biggl\{\frac{1}{\sqrt{2}}z_{\alpha}
\sigma_{j_*^m}+f(0),f(t)\biggr\}, &\quad $\mbox{if $t\geq0$;}$
\vspace*{2pt}\cr
f(t), & \quad$\mbox{otherwise,}$}
\]
and if 
$f(-t_*)\leq-\frac{1}{\sqrt{2}}z_{\alpha}\sigma_{j_*^m}+f(0)$, let
\[
g(t)=\cases{ %
\displaystyle\min\biggl\{-\frac{1}{\sqrt{2}}z_{\alpha}
\sigma_{j_*^m}+f(0),f(t)\biggr\}, &\quad $\mbox{if $t\leq0$;}$
\vspace*{2pt}\cr
f(t), &\quad $\mbox{otherwise.}$}
\]
Then we have
\[
\int_{-1/2}^{1/2}\bigl(f(t)-g(t)
\bigr)^2\,dt\leq\frac{1}{2}z^2_{\alpha}
\frac
{2^{j_*^m-1}}{n}\cdot2^{-j_*^m+2}=\frac{z^2_{\alpha}}{n}.
\]
If $j_*^m=2$, let
\[
g(t)=\cases{ %
\displaystyle\max\biggl\{\frac{1}{\sqrt{2}}z_{\alpha}
\sigma_{j_*^m}+f(0),f(t)\biggr\}, & \quad$\mbox{if $t\geq0$;}$
\vspace*{2pt}\cr
f(t), &\quad $\mbox{otherwise,}$ }
\]
then we have
\[
\int_{-1/2}^{1/2}\bigl(f(t)-g(t)
\bigr)^2\,dt\leq\frac{1}{2}z^2_{\alpha}
\frac
{2^{j_*^m-1}}{n}\cdot\frac{1}{2}\leq\frac{z^2_{\alpha}}{n}.
\]
It then follows that
\[
\omega\biggl(\frac{z_{\alpha}}{\sqrt{n}},f,F_m\biggr)\geq
\frac{1}{\sqrt {2}}z_{\alpha}\sigma_{j_*^m},
\]
and so
\[
L_{\alpha}^*(f, F_m) \geq\biggl(1-\frac{1}{\sqrt{2\pi}z_{\alpha}}\biggr)
\frac
{1}{\sqrt{2}}z_{\alpha}\sigma_{j_*^m}.
\]

\subsection{\texorpdfstring{Proof of Theorem \protect\ref{convex.CI-EL.thm}}{Proof of Theorem 4}}

We shall first prove that the confidence interval $\mathrm{CI}_*^c$ has
guaranteed coverage probability of $1-\alpha$ over $F_c$ and then prove
the upper bound for the expected length.

Note that if $j_*^c >1$, then $E T_{j_*^c -1} \ge\frac{ 2}{3}
z_{\alpha} \sigma_{j_*^c-1}
= \frac{\sqrt2}{3} z_{\alpha} \sigma_{j_*^c}$.
It follows that for $k\ge1$, $ET_{j_*^c -k} \ge2^{k-1/2}\frac{ 1}{3}
z_{\alpha} \sigma_{j_*^c}
= 2^{(3k-1)/2} \frac{1}{3} z_{\alpha} \sigma_{j_*^c - k} $.
Hence
%
\begin{equation}\qquad
\label{too.early} P\bigl(\hat j = j_*^c -k\bigr) \le P(
T_{j_*^c -k} \le z_{\alpha} \sigma_{j_*^c -k}) \le P\biggl(Z \ge
\biggl(\frac{2^{(3k-1)/2}}{3} - 1\biggr) z_{\alpha}\biggr).
\end{equation}

Also for $m \ge0$, $ET_{j_*^c +m} \le2^{-m}\cdot\frac{2}{3}
z_{\alpha} \sigma_{j_*^c}
= 2^{-3m/2} \cdot\frac{2}{3} z_{\alpha} \sigma_{j_*^c +m}
$
and hence
%
\begin{eqnarray}
\label{too.far} P\bigl(\hat j \ge j_*^c + k\bigr) &\le& \prod
_{m=0}^{k-1} P( T_{j_*^c +m} > z_{\alpha}
\sigma_{j_*^c +m})
\nonumber
\\[-8pt]
\\[-8pt]
\nonumber
&\le &\prod_{m=0}^{k-1}
P\biggl( Z > z_{\alpha} \biggl(1 - \frac{2}{3}2^{-3m/2}
\biggr)\biggr).
\end{eqnarray}

To bound the coverage probability note that
%
\begin{eqnarray}
\label{coverage1} P\bigl(f(0) \notin \mathrm{CI}_*^c\bigr) &\le& \sum
_{m \ge3} P\bigl(\hat j = j_*^c -m\bigr) + P\bigl(
\hat j \ge j_*^c + 3\bigr)
\nonumber
\\[-8pt]
\\[-8pt]
\nonumber
&&{}+ \sum_{k = -2}^{2}
P\bigl(f(0) \notin \mathrm{CI}_{j_*^c + k}\bigr).
\end{eqnarray}
It then follows from equation \eqref{too.early} that
\[
P\bigl(\hat j = j_*^c - 3\bigr) \le P\biggl(Z\ge
{13\over3}z_{\alpha}\biggr) \le{7\alpha
\over10\mbox{,}000}
\]
for all $0<\alpha\leq0.2$. It is easy to verify directly that for all
$z\ge1$, $P(Z\ge2z) \le(1/6) P(Z\ge z)$. Furthermore, it is easy to
see that for $k\ge1$, $\frac{2^{(3(k+3)-1)/2}}{3} - 1 \ge2^k
{13\over
3}$ and so
\begin{eqnarray*}
P\bigl(\hat j = j_*^c - 3 - k \bigr) &\le& P\biggl(Z\ge\biggl(
\frac{2^{(3(k+3)-1)/2}}{3} - 1\biggr) z_{\alpha}\biggr) \le P\biggl(Z
\ge2^k {13\over3} z_{\alpha}\biggr)
\\
&\le& 6^{-k} P\biggl(Z\ge{13\over3} z_{\alpha}
\biggr) \le6^{-k}{7\alpha\over10\mbox{,}000}.
\end{eqnarray*}
Hence,
\[
\sum_{m \ge3} P\bigl(\hat j = j_*^c -m
\bigr) = \sum_{k \ge0} P\bigl(\hat j = j_*^c
- 3 - k \bigr) \le{7\alpha\over10\mbox{,}000 }\sum_{k \ge0}6^{-k}
\le{7\alpha
\over
5000 }.
\]

Note that \eqref{too.far} yields that
\[
P\bigl(\hat j \ge j_*^c + 3\bigr) \le P\biggl(Z\ge
{1\over3}z_{\alpha}\biggr) \cdot P\biggl(Z\ge \biggl(1-
{1\over3\sqrt{2}}\biggr)z_{\alpha}\biggr) \cdot P\biggl(Z\ge
{11\over12}z_{\alpha}\biggr) \le {\alpha\over6.4}
\]
for all $0<\alpha\le0.3$. It now follows from \eqref{coverage1} that
\[
P\bigl(f(0)\in \mathrm{CI}_*^c\bigr) = 1 - P\bigl(f(0)\notin
\mathrm{CI}_*^c\bigr) \ge1 - \biggl({7\alpha\over
5000 }+
{\alpha\over6.4} + 5\times{\alpha\over6}\biggr)\ge1- \alpha.
\]

We now turn to the upper bound for the expected length.
Note that
%
\begin{equation}
E_f L\bigl(\mathrm{CI}_*^c\bigr) \le\sum
_{j=1}^{\infty} \bigl( z_{\alpha}+ (\sqrt{5} +
\sqrt{2}) z_{\alpha/12}\bigr) \sigma_j \cdot P(\hat j =j).
\end{equation}
Hence
%
\begin{eqnarray}
\label{EL} E_f L\bigl(\mathrm{CI}_*^c\bigr) &\le& \bigl(
z_{\alpha}+ (\sqrt{5} + \sqrt{2}) z_{\alpha/12}\bigr) \sigma
_{j_*^c}
\nonumber
\\[-8pt]
\\[-8pt]
\nonumber
&&{}\times\Biggl(P\bigl( \hat j \le j_*^c\bigr) + \sum
_{k=1}^{\infty} 2^{k/2} P\bigl(\hat j
=j_*^c +k\bigr) \Biggr).
\end{eqnarray}
Set $w_k = 2^{k/2} - 2^{(k-1)/2}$ for $k\geq1$. Then it is easy to see that
\[
S = P\bigl( \hat j \le j_*^c\bigr) + \sum
_{k=1}^{\infty}2^{k/2} P\bigl(\hat j
=j_*^c +k\bigr) = 1 + \sum_{k=1}^{\infty}w_k
P\bigl(\hat j \ge j_*^c +k\bigr).
\]
It then follows from \eqref{too.far} that
\[
S \le1 + \sum_{k=1}^{\infty}w_k
\prod_{m=0}^{k-1} P\biggl( Z >
z_{\alpha} \biggl(1 - \frac{2}{3}\frac{1}{2^{3m/2}}\biggr)\biggr).
\]
The right-hand side is clearly increasing in $\alpha$. Direct numerical
calculations show that for $\alpha=0.2$,
\[
\sum_{k=1}^{\infty}w_k \prod
_{m=0}^{k-1} P\biggl( Z > z_\alpha
\biggl(1 - \frac
{2}{3}\frac{1}{2^{3m/2}}\biggr)\biggr) \le0.25.
\]
It then follows directly from \eqref{EL} that
\[
EL\bigl(\mathrm{CI}_*^c\bigr) \le1.25 \bigl(z_\alpha+ (\sqrt{5} +
\sqrt{2}) z_{\alpha/
12}\bigr) \sigma_{j_*^c}. 
\]

\subsection{\texorpdfstring{Proof of Theorem \protect\ref{EL.bound2.thm}}{Proof of Theorem 5}}

Note that if $j_*^c >1$,
then $E T_{j_*^c -1} \ge\frac{ 2}{3} z_{\alpha} \sigma_{j_*^c-1}
= \frac{\sqrt2}{3} z_{\alpha} \sigma_{j_*^c}$,
and hence there is a $t_*$ satisfying $0<t_* \le2^{-j_*^c +1 }$ such that
$f_s(t_*) = \frac{\sqrt2}{3} z_{\alpha} \sigma_{j_*^c}$, where
$f_s(t)=\frac{f(t)+f(-t)}{2}-f(0)$.
Let $g$ be defined by
\[
g(t) = f(t) 1\bigl(|t| > t_*\bigr) + \biggl(f_s(t_*) + \frac{f(t_*) - f(-t_*)}{2t_*} t
\biggr)1\bigl(|t| \le t_*\bigr).
\]
There is also a $g$ as in the proof of Lemma 5 in our other paper with
$g(0) = f_s(t_*)$
for which
\[
\int_{-1/2}^{1/2} \bigl(g(t) - f(t)
\bigr)^2 \,dt \le\frac{9}{4} f^2_s(t_*)
t_* 
\leq\frac{z_{\alpha}^2}{n}.
\]
If $j_*^c=1$, then let $g(t)=f(t)+\frac{\sqrt2}{3} z_{\alpha} \sigma
_{j_*^c}$, and then we have
\[
\int_{-1/2}^{1/2} \bigl(g(t) - f(t)
\bigr)^2 \,dt \le\frac{2}{9}z_{\alpha
}^2
\sigma_1^2 \leq\frac{z_{\alpha}^2}{n}.
\]
It then follows that
\[
\omega\biggl( \frac{z_{\alpha}}{\sqrt{n}},f,F_c\biggr) \ge
\frac{ \sqrt
2}{3}z_{\alpha} \sigma_{j_*^c},
\]
and so
\[
L_{\alpha}^{*}(f, F_c) \ge\biggl(1 -
\frac{1}{\sqrt{2\pi} z_{\alpha}}\biggr) \frac
{\sqrt2}{3}z_{\alpha} \sigma_{j_*^c}.
\]

%



\printaddresses

\end{document}